\documentclass[12pt,reqno,oneside]{amsart}
 \usepackage{verbatim,amsmath,amssymb,cite,epsfig,xspace}
\usepackage{color,cite,graphicx}
   
\numberwithin{equation}{section}

\theoremstyle{plain}
\newtheorem{theorem}{Theorem}[section]
\newtheorem{lemma}{Lemma}[section]

\theoremstyle{definition}
\newtheorem{definition}{Definition}[section]

\theoremstyle{remark}

\begin{document}

\title[]
{Liouville theorems in unbounded domains for the time-dependent Stokes System }
\author{Hao Jia, Gregory Seregin, Vladim\'{i}r \v{S}ver\'{a}k}
\address{}
\email{}
\maketitle

\centerline{\sl Dedicated to Professor Peter Constantin on the occasion of his 60th birthday. }

\begin{abstract}{In this paper, we characterize bounded ancient solutions to the time-dependent Stokes system with zero boundary value in various domains, including the half-space.}
\end{abstract}

\begin{section}{Introduction}
 In this paper, we show that any solution $u(x,t)\in L^{\infty}(\Omega\times (-\infty,0))$ to
\begin{eqnarray}
\left.\begin{array}{rl}
        \partial_t u-\Delta u+\nabla p &=0\\
         \mbox{div}~~u&=0
       \end{array}\right\}&& \quad \mbox{in $\Omega\times (-\infty,0)$, and}\label{eq:main1}\\
        u|_{\partial \Omega}=0&& \label{eq:main2}
\end{eqnarray}
 satisfies
\begin{eqnarray}\label{eq:result}
 u(x,t)=\left\{\begin{array}{ll}
                0 &~~~ \mbox{if $\Omega\subset R^n$ is a bounded domain and $n\ge 2$,}\\
            u(t) &~~~ \mbox{if $\Omega=R^n$ and $n\ge 2$,}\\
    u(t,x_n) ~~~\mbox{with}~~~u_n=0 &~~~ \mbox{if $\Omega$ is half space $x_n>0$ and $n\ge 2$,}\\
      a(t)+O(\frac{1}{|x|^{n-2}}) ~~~\mbox{as}~~~|x|\to \infty &~~~\mbox{if $\Omega\subset R^n$ is an exterior domain with $n\ge 3$. }
               \end{array}\right.
\end{eqnarray}

 Throughout the paper we assume that the domains $\Omega$ have smooth boundary. One has to be somewhat careful with the definition of the boundary condition $u|_{\partial\Omega}=0$ since a-priori we only assume $u$ to be bounded, with no further regularity assumptions. The usual definition is the following:
\begin{definition}
 We call $u(x,t)\in L^{\infty}(\Omega\times (0,\infty))$ a very weak ancient solution to equations (\ref{eq:main1})(\ref{eq:main2}) if
\begin{equation}\label{eq:main3}
 \int_{-\infty}^0\int_{\Omega}u(x,t)(\partial_t \phi+\Delta \phi)(x,t)dxdt=0
\end{equation}
for any $\phi\in C_c^{\infty}(\overline{\Omega}\times (-\infty,0))$ with $\mbox{div}~~\phi=0$, $\phi|_{\partial \Omega\times (-\infty,0)}=0$; and
\begin{equation}\label{eq:main4}
 \int_{-\infty}^0\int_{\Omega}u(x,t)\nabla \psi(x,t)dxdt=0
\end{equation}
for any $\psi\in C_c^{\infty}(\overline{\Omega}\times (-\infty,0))$.
\end{definition}
Solutions defined in this way  are often called very weak solutions in the literature and we also use this terminology. For smooth solutions the definition coincides with the usual one, as one can easily check by integration by parts. Our main result is as follows:
\begin{theorem}
 Let $u$ be a bounded very weak ancient solution to equations (\ref{eq:main1}),(\ref{eq:main2}), in the sense that $u$ satisfies equations (\ref{eq:main3}),(\ref{eq:main4}). Then $u$ is given by (\ref{eq:result}).
\end{theorem}

\noindent
\textbf{Remarks}: The results are essentially sharp. This is obvious in the cases when $\Omega$ is $R^n$ or a bounded domain. In the case of a half space, one can take $u=(u_1(t,x_n),\dots,u_{n-1}(t,x_n),0)$, $u_i$ verifies: $\partial_tu_i-\partial_n^2u_i=f_i(t), ~~u_i|_{x_n=0}=0$ for $1\leq i\leq n-1$, where $f_i\in C_c^{\infty}(-\infty,0)$. Then $u$ is a solution to equations (\ref{eq:main1}),(\ref{eq:main2}). This is the example given in \cite{SV10}. In exterior domain, the decay rate we obtain is as good as that of the fundamental solution of steady Stokes equations. \\

Our work is motivated by boundary regularity for the Navier-Stokes equations. The main interest is in the case of the half-space, the other case are included for completeness. The connection between regularity and Liouville-type theorems is of course classical. In the context of the Navier-Stokes equations it is discussed for example in \cite{KNSS, SS}. In a recent note \cite{SV10} a bounded shear flow for unsteady Stokes equations is constructed which is not fully regular although the boundary value is zero. This example simplifies earlier constructions of \cite{kK01}. The lack of boundary regularity in the time-dependent case is in contrast with the case of steady Stokes equations, see e.g. \cite{Kang}. In the time-dependent Stokes equations and Navier-Stokes equations, one usually treats pressure as an auxiliary variable, determined by $u$. Such treatment is valid as long as we have some decay of $u$ at spatial infinity. On the other hand, it has been known that in unbounded domains, if we do not assume decay of $u$, the pressure may act as an external force `driving' the fluid motion, as in the case of \cite{SV10}. In such situations, we lose boundary regularity even with the vanishing boundary conditions. In this context, our result could be understood as showing that the solutions in \cite{SV10} are in some sense the only obstacle to full boundary regularity (in suitable solution classes).

 The paper is organized as follows. In section 2 we introduce some technical lemmas to be used below. Section 3 deals with the simple cases when $\Omega$ is a bounded domain or the whole space. Section 4 and 5 deal with the more subtle cases when $\Omega$ is a half space or an exterior domain. For the exterior domains, we use a standard extension argument together with some estimates of linear Stokes system. For the half space, which is the most interesting case,  we use Fourier transform. There is also a proof based on duality arguments, which requires some additional point-wise estimates of solutions to linear Stokes system in half space. The estimates may be of independent interest, but the calculations are somewhat lengthy. This alternative proof will appear elsewhere.\\

\noindent
\textbf{Notation}. We will use standard notations. For example, $\Omega$ will be one of the four types of domain in $R^n$ mentioned above. $B_r(x_0):=\{x\in R^n||x-x_0|<r\}$, $Q_r(x_0,t_0):=B_r(x_0)\times (t_0-r^2,t_0)$, $Q_r:=Q_r(0,0)$. $C$ denotes an absolute positive number. $C(\alpha, \lambda,\cdots)$ denote a positive number depending only on $\alpha,\lambda,\cdots$, $\overline{A}$ denotes the closure of $A$, $A\Subset \mathcal{O}$ means the closure of $A$ is a compact subset of $\mathcal{O}$, $\partial_i:=\frac{\partial}{\partial x_i}$.

\end{section}

\begin{section}{Some technical lemmas}
 In the sequel we will make use standard mollifications. For completeness we include the following standard lemma:
\begin{lemma}{}\label{lm:mollification}
 Let $\Omega$ be as above, $u\in L^{\infty}_{x,t}(\Omega\times(-\infty,0))$. Take a standard smooth cutoff function $\eta(t)$ with $supp~~ \eta\Subset(0,1)$ and $\int\eta=1$. For each $\epsilon>0$, we define $u^{\epsilon}$ as a distribution in $\Omega\times (-\infty,0)$ in the following way,
\begin{equation}
 (u^{\epsilon},\phi)=\int_{-\infty}^0\int_{\Omega}u(x,t)\int_{-\infty}^0\frac{1}{\epsilon}\eta(\frac{s-t}{\epsilon})\phi(x,s)dsdxdt
\end{equation}
for any smooth $\phi$ with $supp~~\phi\Subset\Omega\times(-\infty,0)$.\\
Then $u^{\epsilon}$ is a bounded function with bounded distributional derivatives $\partial_t^ku^{\epsilon}$, $k=0,1,2,\cdots$. Moreover, we have the following estimates:
\begin{equation}
 \|\partial_t^ku^{\epsilon}\|_{L^{\infty}(\Omega\times(-\infty,0))}\leq C\epsilon^{-k}\|u\|_{L^{\infty}(\Omega\times(-\infty,0))}\mbox{.}
\end{equation}

\end{lemma}

\noindent
\textbf{Proof and Remarks}: The proof follows immeditely from well-known properties of convolution.  We note that due to our special choice of the support of $\eta$, the mollified function $u^{\epsilon}$ is still defined in $\Omega\times(-\infty,0)$. It is clear from definition that $u^{\epsilon}$ converge weakly$^\ast$ to $u$ in $L^{\infty}(\Omega\times(-\infty,0))$. It is also clear that, after possibly changing the value of $u^{\epsilon}$ on a set of measure zero, the map $t\to u^{\epsilon}(\cdot,t)$ is continuous from $(-\infty,0)$ to $L^p(K)$ for any $K\Subset\Omega$, $1<p<\infty$. \\

Let $u$ be a bounded distributional solution to the linear Stokes equations (\ref{eq:main1}) in $Q_1$ with some distribution $p$. It is well known that we have regularity of $u$ in $x$, for almost every $t$ in $Q_{1/2}$. We can not, however, expect to have any regularity in $t$ for $u$, or any reasonable estimate on $p$ in general, assuming only that $u$ is bounded in $Q_1$. This point is usually illustrated with the example where $u(t,x)=f(t)$, $p=-f'(t)\cdot x$. Here $f(t)$ is bounded, but $f'(t)$ can be arbitrarily large. On the other hand, if we assume some estimate on $\partial_tu$, then we can improve estimates on $p$. The following lemma summarizes the above discussion.
\begin{lemma}\label{lm:linearregularity}
 Let $u$ be a bounded distributional solution to linear unsteady Stokes equations in $Q_1$. Let $\|u\|_{L^{\infty}_{x,t}(Q_1)}\leq 1$. Then for any multi-index $\alpha$ with $|\alpha|\ge 0$, $\|\partial_x^{\alpha}u\|_{L^{\infty}_{x,t}(Q_{1/2})}\leq C(n,\alpha)$. If in addition $\|\partial_tu\|_{L^{\infty}_{x,t}(Q_1)}\leq M$, then there exists a pressure field $p(x,t)$ such that (\ref{eq:main1}) is satisfied and $\|\partial_x^{\alpha}p\|_{L^{\infty}_{x,t}(Q_{1/2})}\leq C(\alpha,M,n)$ for any multi-index $\alpha$ with $|\alpha|\ge 0$.
\end{lemma}

\noindent
\textbf{Proof}: For the first part of the lemma, note that the vorticity $\omega_{ij}:=\partial_iu_j-\partial_ju_i$, $1\leq i,j\leq n$, satisfies heat equation $\partial_t\omega_{ij}-\Delta \omega_{ij}=0$ in $Q_1$. Thus $\omega_{ij}$ is smooth with all derivatives bounded by constants depending only on $n$ in $Q_{3/4}$. From the divergence free condition, we get $\Delta u_i=-\sum_{j=1}^n\partial_j\omega_{ij}$. Then the first part of the lemma follows from interior estimate of Laplace equations. For the second part, note that $\|\nabla p\|_{L^{\infty}_{x,t}(Q_{3/4})}\leq C(n,M)$ from the assumption on $\partial_tu$ and first part of the lemma. Since we also have $\Delta p=0$, the estimate follows. \\

\smallskip
\noindent
\textbf{Remark}: The pressure is only determined up to an arbitrary function  of $t$. (If we change $p$ to $p+c(t)$, equation (\ref{eq:main1}) is not affected.) In estimates below we will usually  assume a suitable choice of $c(t)$.\\

We shall need the following extension result (which is interesting in its own right) below.
\begin{lemma}{(Extension of divergence-free vector field)}\label{lm:extension}
 For any smooth compactly supported vector field $g=(g_1,\cdots,g_{n-1},0)$ in $R^{n-1}$, there exists a smooth divergence free vector field $\phi=(\phi_1,\cdots,\phi_n)$ with ${\rm supp~~}\phi\Subset \overline{R^n_{+}}$ such that $\phi|_{x_n=0}=0$ and $\frac{\partial \phi}{\partial x_n}|_{x_n=0}=g$.
\end{lemma}

\noindent
\textbf{Proof}: We seek $\phi$ in the form of $\phi_i=\sum_{j=1}^n\partial_j w_{ij}$, with some $w_{ij}\in C^{\infty}_c(\overline{R^n_{+}})$ and $w_{ij}=-w_{ji}$ for $1\leq i,j\leq n$. Note that under such conditions on $w_{ij}$, $\mbox{div~~}\phi=0$ is automatically satisfied. To satisfy boundary conditions for $\phi$, we need:
\begin{equation} \label{eq:boundary1}
 \sum_{j=1}^n\frac{\partial w_{ij}}{\partial x_j}|_{x_n=0}=0\quad\mbox{for $1\leq i\leq n$,}
\end{equation}
and
\begin{equation} \label{eq:boundary2}
 \sum_{j=1}^n\frac{\partial^2w_{ij}}{\partial x_n\partial x_j}|_{x_n=0}=g_i\quad \mbox{for $1\leq i\leq n$, $g_n$=0.}
\end{equation}
It is easy to verify that the $n$-th equation in (\ref{eq:boundary2}) is also automatically satisfied once the rest of the equations in the above are satisfied. To satisfy equations (\ref{eq:boundary1})(\ref{eq:boundary2}), we first require $w_{ij}|_{x_n=0}=0$, for all $1\leq i,j\leq n$. Then equations (\ref{eq:boundary1}) become $\frac{\partial w_{in}}{\partial x_n}|_{x_n=0}=0$ for $1\leq i\leq n$. We further require that $w_{ij}=0$ if $1\leq i,j\leq n-1$, then equations (\ref{eq:boundary2}) reduce to $\frac{\partial^2w_{in}}{\partial^2x_n}|_{x_n=0}=g_i$ for $1\leq i\leq n-1$. Summarizing the above analysis, it is sufficient to find $w_{in}\in C_c^{\infty}(\overline{R^n_{+}})$ for $1\leq i\leq n-1$, $w_{in}=-w_{ni}$ satisfying
\begin{equation}
 \left.\begin{array}{rl}
       w_{in}|_{x_n=0}&=0\\
    \frac{\partial w_{in}}{\partial x_n}|_{x_n=0}&=0\\
\frac{\partial^2w_{in}}{\partial^2x_n}|_{x_n=0}&=g_i
      \end{array}\right\}\quad\mbox{for $1\leq i\leq n-1$.}
\end{equation}
It is clear that we can always find such $w_{in}$. Thus $\phi$ satisfying conditions in the lemma exists.\\

We collect some facts about the operator $|\nabla|$ which will be used in our proofs. \\

For $f\in \mathcal{S}(R^n)$, we define $|\nabla|f(x)=(|\xi|\hat{f}(\xi))^{\vee}(x)$ where we have used Fourier transform $\hat{f}(\xi)=\frac{1}{(2\pi)^{n/2}}\int_{R^n}e^{-ix\cdot \xi}f(x)dx$ and inverse Fourier transform $\check{f}(x)=\frac{1}{(2\pi)^{n/2}}\int_{R^n}e^{ix\cdot \xi}f(\xi)d\xi$. One can write $|\nabla|f=\sum_{j=1}^n-\frac{\partial_j}{|\nabla|}\partial_jf=\sum_{j=1}^n R_j\partial_j f$, where $R_j$ denotes the Riesz transform. Clearly $|\nabla|$ can be considered as a continuous operator from $\mathcal{S}(R^n)\to L^1(R^n)$. By duality, we can extend $|\nabla|$ to act on $L^{\infty}(R^n)$ according to the usual fromula
$\langle|\nabla|f,\phi\rangle=\langle f,|\nabla|\phi\rangle$ for any  $\phi\in \mathcal{S}(R^n)$.\\

 We recall the following obvious continuity result.
\begin{lemma}\label{L:lemmano1}
 Let $|\nabla|$: $L^{\infty}(R^n)\longmapsto \mathcal{S}'(R^n) $ be defined as above. If $u_m\in L^{\infty}$, with $u_m$ converges weakly$^\ast$ to $u$ in
$L^{\infty}$ (viewed as the dual of $L^1(R^n)$), then $|\nabla|u_m$ converges to $|\nabla|u$ in $\mathcal{S}'(R^n) $
\end{lemma}

\noindent
\textbf{Proof}: This follows directly from the definitions.\\

Recall the definition of H\"{o}lder norm in $R^n$: $\|u\|_{C^{m,\alpha}(R^n)}:=\sum_{i=0}^m\sup_{|\beta|\leq m}\sup_{x\in R^n}|\partial^{\beta}u(x)|+\sup_{|\beta|=m}\sup_{x,y\in R^n,~~x\neq y}\frac{|\partial^{\beta}u(x)-\partial^{\beta}u(y)|}{|x-y|^{\alpha}}$  for any $m\ge 0$, $0<\alpha<1$. The H\"{o}lder space $C^{m,\alpha}(R^n)$ is consisted of all $u$ with $\|u\|_{C^{m,\alpha}}<\infty$. We will use  the following estimate:
\begin{lemma}\label{L:lemmano2}
 $|\nabla|:C^{m+1,\alpha}(R^n)\mapsto C^{m,\alpha}(R^n)$ is bounded, for $m\geq 1,~~0<\alpha<1$.
\end{lemma}

\noindent
\textbf{Proof}: This follows from the representation $|\nabla|=\sum_j R_j\partial_j$ and the Schauder estimates for the Riesz transform.

\medskip
\noindent
We will denote by $|\nabla'|$  the analogue of $|\nabla|$ acting only on the variables $x_1,\dots,x_{n-1}$.
For Schwartz function $f$, $|\nabla'|f(x',x_n)=(2\pi |\xi'|\hat{f}(\xi',x_n))^{\vee}(x')$, where the Fourier transform and inverse Fourier  transform are both with respect to the first $n-1$ variables. From definition, it is clear if $f(x',x_n,t)\in L^{\infty}(R^{n-1}\times(x_1,x_2)\times(t_1,t_2))$, then $|\nabla'|f\in \mathcal{D}'(R^{n-1}\times(x_1,x_2)\times(t_1,t_2))$. Moreover, if $\partial_{x'}^l\partial_n^k\partial_t^mf\in L^{\infty}$, then $|\nabla'|\partial_{x'}^l\partial_n^k\partial_t^mf=\partial_{x'}^l\partial_n^k\partial_t^m|\nabla'|f$ in $\mathcal{D}'(R^{n-1}\times(x_1,x_2)\times(t_1,t_2))$.\\

For bounded harmonic functions in the upper half spaces, we have the following result (see also \cite{sU87}, for example).
\begin{lemma}\label{lm:lemmano3}
 Let $f$ be a bounded harmonic function in the upper half space $R^n_{+}$, we have $(\partial_n f+|\nabla'|f)(x)=0$ in the sense of distributions in $R^n_{+}$.
\end{lemma}

\noindent
\textbf{Remarks:} By classical regularity for harmonic functions and lemma \ref{L:lemmano2}, we see both $\partial_nf$ and $|\nabla'|f$ are smooth functions in the interior of $R^n_+$.

\medskip
\noindent
\textbf{Proof}: Let $P(x,y)$ be the Poisson kernel. By classical representation results there exists a $g\in L^{\infty}(R^{n-1})$,
such that $f(x)=\int_{R^{n-1}} P(x,y)g(y)\,dy.$ By approximation and continuity properties of $|\nabla'|$ we can assume without loss of generality that $g$ is smooth and compactly supported. 
 Applying Fourier transform in the $x_1,\dots,x_{n-1}$ variables, we have $\hat f(\xi',x_n)=\hat g(\xi')e^{-|\xi'|x_n}$ and the result follows.
\end{section}

\begin{section}{The cases $\Omega=R^n$ or a bounded domain}
In this section, we first deal with the (easy) cases when $\Omega=R^n$ or $\Omega$ is a bounded domain. Recall that our goal is to show that bounded very weak ancient solutions to (\ref{eq:main1})(\ref{eq:main2}) are given by (\ref{eq:result}).\\

\medskip
\noindent
1. $\Omega=R^n$.\\
In this case, it is not difficult to see that equations (\ref{eq:main3})(\ref{eq:main4}) are equivalent to
\begin{equation*}
 \left.\begin{array}{rl}
       \partial_tu-\Delta u+\nabla q&=0\\
       \mbox{div}~~u&=0
      \end{array}\right\}\quad\mbox{in $R^n\times(-\infty,0)$}
\end{equation*}
in the sense of distributions for some $q\in \mathcal{D}'(R^n\times(-\infty,0))$.\\
For $1\leq i,j\leq n$, let $\omega_{ij}=\partial_ju_i-\partial_iu_j$, then clearly $\partial_t\omega_{ij}-\Delta \omega_{ij}=0$ in $R^n\times(-\infty,0)$. Since $\omega_{ij}$ are bounded in some negative Sobolev space, we immediatly get $\omega_{ij}$ are bounded functions from parabolic regularity. Thus $\omega_{ij}$ are so called bounded ancient solution to heat equation, and consequently $\omega_{ij}=$constants $c_{ij}$.  Since $u$ is divergence free, we get $\Delta u_i=-\sum_{j=1}^n\partial_j\omega_{ij}=0$ in $R^n\times(-\infty,0)$. Therefore $u(t,x)=f(t)$ for some bounded measurable $f$ a.e t. This completes the proof when $\Omega=R^n$.\\

\medskip
\noindent
2. $\Omega$ is a bounded domain.\\
In this case our goal is to show that bounded very weak ancient solutions $u$ to (\ref{eq:main1})(\ref{eq:main2}) are identically 0. We use a duality argument as follows. For any $f\in C_c^{\infty}(\Omega\times(0,+\infty))$, let $\tilde{\phi}$ solve
\begin{eqnarray*}
 \left.\begin{array}{rl}
       \partial_t\tilde{\phi}-\Delta \tilde{\phi}+\nabla q&=f\\
       \mbox{div}~~\tilde{\phi}&=0
      \end{array}\right\}&&\quad\mbox{in $\Omega\times(0,\infty)$,}\\
\tilde{\phi}(\cdot,t)|_{\partial\Omega}=0\mbox{.}&&
\end{eqnarray*}
The existence, uniqueness and regularity of such solutions are well known, one can see e.g \cite{Galdi}. Moreover, we have $\lim_{t\to \infty}\|\tilde{\phi}(\cdot,t)\|_{L^2(\Omega)}=0$ (the decay is actually exponential). Take a standard smooth cutoff function $\eta(t)$ with $\eta(t)=0$ for $t>2$. For any $R>0$, let $\phi_R(x,t)=\eta(-\frac{t}{R})\tilde{\phi}(x,-t)$ for $t\in (-\infty,0)$. Then from equations (\ref{eq:main3})(\ref{eq:main4}) we obtain
\begin{eqnarray*}
 &&0=\int_{-\infty}^0\int_{\Omega}u(x,t)(\partial_t\phi_R+\Delta\phi_R)dxdt\\
&&=\int_{-\infty}^0\int_{\Omega}u(x,t)(-\partial_t\tilde{\phi}+\Delta\tilde{\phi})(x,-t)\eta(-\frac{t}{R})dxdt\\
&&-\frac{1}{R}\int_{-\infty}^0\int_{\Omega}u(x,t)\eta'(-\frac{t}{R})\tilde{\phi}(x,-t)dxdt\\
&&=-\int_{-\infty}^0\int_{\Omega}u(x,t)f(x,-t)\eta(-\frac{t}{R})dxdt-\frac{1}{R}\int_{-\infty}^0\int_{\Omega}u(x,t)\eta'(-\frac{t}{R})\tilde{\phi}(x,-t)dxdt\\
&&+\int_{-\infty}^0\int_{\Omega}u(x,t)\eta(\frac{t}{R})\nabla q(x)dxdt
\end{eqnarray*}
Using the fact that $f$ is compactly supported in $t$, $q$ is smooth in $x$, $u$ is bounded and $\lim_{t\to \infty}\|\tilde{\phi}(\cdot,t)\|_{L^1(\Omega)}=0$ (since $\Omega$ is bounded), we can send $R\to \infty$ and obtain
\begin{equation*}
 \int_{-\infty}^0\int_{\Omega}u(x,t)f(x,-t)dxdt=0\mbox{.}
\end{equation*}
Since $f$ is arbitrary, we must have $u\equiv 0$.
\end{section}
\begin{section}{The case $\Omega=R^{n}_{+}$}
 Now let us deal with the more subtle case when $\Omega$ is a half space. In fact one can still use the idea of duality as in the case of bounded domains. In this case, however, one has to study the decay property of solution to the linear Stokes equations quite carefully. One also has to appropriately localize $\tilde{\phi}$ (assuming notations from the last section) since in (\ref{eq:main3}) the test function $\phi$ is required to be of compact support. The authors have obtained a proof using such a method, which will appear elsewhere.\\
Here we take a different approach based on the Fourier transform in which the calculations are simpler.\\

Let $u$ be as above, take a smooth mollifier $\eta(x',t)$ with $supp~~\eta\Subset B_1(0)\times(0,1)\subseteq R^{n-1}_{x'}\times R_t$ and $\int \eta=1$. We define the mollified vector field $u^{\epsilon}$ similar as before, again by duality: for any smooth $\phi$ with $supp~~\phi\Subset R^n_{+}\times(-\infty,0)$,
\begin{equation*} (u^{\epsilon},\phi):=\int_{-\infty}^0\int_{R^n_{+}}u(x,t)\int_{R^{n-1}}\int_{-\infty}^0\epsilon^{-n}\eta(\frac{y'-x'}{\epsilon},\frac{s-t}{\epsilon})\phi(y',x_n,s)dy'dsdxdt\mbox{.}
\end{equation*}
Again similar as before, one can show $u^{\epsilon}$ is bounded with bounded distributional derivatives $|\partial_t^k\partial_{x'}^{\alpha}u^{\epsilon}|\leq C(k,\alpha,n)\epsilon^{-k-|\alpha|}\|u\|_{L^{\infty}_{x,t}}$. We have the following result:
\begin{lemma}
 Let $u$ be a bounded very weak ancient solution to (\ref{eq:main1})(\ref{eq:main2}) in $R^n_{+}\times(-\infty,0)$, $u^{\epsilon}$ is defined as above. Then $u^{\epsilon}$ is smooth with all derivatives bounded in $\overline{R^n_{+}}\times (-\infty,0)$. Moreover, $u^{\epsilon}$ still satisfies equations (\ref{eq:main3}),(\ref{eq:main4}) and $u^{\epsilon}(\cdot,t)|_{x_n=0}=0$.
\end{lemma}

\noindent
\textbf{Proof}: From equations (\ref{eq:main3})(\ref{eq:main4}) and definition of $u^{\epsilon}$ we see
\begin{equation}\label{eq:no1}
 \int_{-\infty}^0\int_{\Omega}u^{\epsilon}(x,t)(\partial_t \phi+\Delta \phi)(x,t)dxdt=0
\end{equation}
for any $\phi\in C_c^{\infty}(\overline{\Omega}\times (-\infty,0))$ with $\mbox{div}~~\phi=0$, $\phi|_{\partial \Omega\times (-\infty,0)}=0$; and
\begin{equation}\label{eq:no2}
 \int_{-\infty}^0\int_{\Omega}u^{\epsilon}(x,t)\nabla \psi(x,t)dxdt=0
\end{equation}
for any $\psi\in C_c^{\infty}(\overline{\Omega}\times (-\infty,0))$. These clearly imply
\begin{equation}
 \left.\begin{array}{rl}
        \partial_tu^{\epsilon}-\Delta u^{\epsilon}+\nabla\cdot q&=0\label{eq:mollified}\\
        \mbox{div}~~u^{\epsilon}&=0
       \end{array}\right\} \quad\mbox{in $\mathcal{D}'(R^n_{+}\times(-\infty,0))$.}
\end{equation}
We first show $u^{\epsilon}$ is smooth up to boundary $\{x_n=0\}$. From the differentiablity property of $u^{\epsilon}$ in $x',t$, we see $q$ is well defined for each $t\in(\infty,0)$ modulo some $c(t)$. Moreover, from $\Delta q=0$ and elliptic estimates, we know $q$ is smooth in $x$ away from boundary $\{x_n=0\}$. Now let us rewrite the $n$-th equation of (\ref{eq:mollified}) as
\begin{equation*}
 0=\partial_tu^{\epsilon}_n-\Delta u^{\epsilon}_n+\partial_nq=\frac{\partial}{\partial x_n}(q-\partial_nu^{\epsilon}_n)-\Delta_{x'}u^{\epsilon}_n+\partial_tu^{\epsilon}_n\mbox{.}
\end{equation*}
Note that $\partial_n u^{\epsilon}_n=-\sum_{i=1}^{n-1}\partial_i u^{\epsilon}_i$ is bounded up to $x_n=0$. We see $q-\partial_nu^{\epsilon}_n$ is bounded up to boundary, thus $q$ is bounded up to boundary $\{x_n=0\}$. The same argument also shows $\nabla_{x'}^{\alpha}q$ is bounded up to boundary. Use $\Delta q=0$ we obtain that $q$ is smooth in spatial variables up to $x_n=0$. Then from $\partial_n^2u^{\epsilon}=\partial_tu^{\epsilon}-\Delta_{x'}u^{\epsilon}+\nabla q$ we get $u\in C^2$. By differentiating the equations in $x_n$ and applying similar arguments we obtain smoothness of $u^{\epsilon}$. Next we show $u^{\epsilon}|_{x_n=0}=0$. Since $u^{\epsilon}$ is smooth in $\overline{R^n_{+}}$, we can use equations (\ref{eq:main3})(\ref{eq:main4}) and integration by parts in equations (\ref{eq:no1})(\ref{eq:no2}) to obtain:
\begin{eqnarray}
 \int_{-\infty}^0\int_{R^{n-1}}u^{\epsilon}_n\psi dxdt&=&0\mbox{,}\\
\int_{-\infty}^0\int_{R^{n-1}}u^{\epsilon}\phi dxdt&=&0\mbox{.}
\end{eqnarray}
Clearly $\psi$ can be arbitrary smooth compactly supported function, thus $u^{\epsilon}_n|_{x_n=0}\equiv 0$. By lemma \ref{lm:extension}, $\phi|_{x_n=0}$ can be any smooth compactly supported vector field with zero $n$-th component, thus $u^{\epsilon}_i|_{x_n=0}\equiv 0$ for $1\leq i\leq n-1$. Therefore, $u^{\epsilon}|_{x_n}\equiv 0$.\\

Now we can prove our main theorem in this section.
\begin{theorem}
 Let $u$ be a bounded very weak ancient solution to equations (\ref{eq:main1})(\ref{eq:main2}) in $R^n_{+}\times (-\infty,0)$, then we must have $u(x,t)=u(x_n,t)$ and $u_n\equiv 0$.
\end{theorem}

\noindent
\textbf{Proof}: By the above results, it is clear we only need to prove our theorem in the case $u(t,x)$
is smooth up to boundary, with all derivatives bounded and $u|_{x_n=0}=0$. Then we see $\frac{\partial p}{\partial x_i}$ is bounded, for $1\leq i\leq n$. Since
$\Delta \frac{\partial p}{\partial x_i}=0$, from lemma \ref{lm:lemmano3} we get $(\partial_n+|\nabla'|)\frac{\partial p}{\partial x_i}=0$.
Applying operator $(\partial_n+|\nabla'|)$ to the $n$-th equation of (\ref{eq:main1}), noting the commutativity of various Fourier multipliers (since $u$ is smooth), we infer that $(\partial_n+|\nabla'|)u_n$
satisfies the heat equation in $R^n\times (-\infty,0)$ with $(\partial_n+|\nabla'|)u_n$ bounded and $(\partial_n+|\nabla'|)u_n|_{x_n=0}=0$ (since $\partial_nu_n=-\sum_{i=1}^{n-1}\partial_iu_i$ and $u|_{x_n=0}=0$). Thus by
Liouville's theorem for heat equation in a half space, we get $(\partial_n+|\nabla'|)u_n\equiv 0$, and consequently $\Delta u_n=0$. Since we also have
$u_n|_{x_n=0}=0$, we see $u_n\equiv 0$. Therefore $\frac{\partial p}{\partial x_n}=0$. Thus $|\nabla'|\frac{\partial p}{\partial x_i}=0$
for $1\leq i \leq n$. Again applying operator $|\nabla'|$ to the first $n-1$ equations of (\ref{eq:main1}), we get $|\nabla'|u'$ satisfies heat equation in $R^n\times (-\infty,0)$ with $(|\nabla'|u')|_{x_n=0}=0$ and
$|\nabla'|u'$ bounded. Using Liouville's theorem for heat equation in a half space again, we obtain $|\nabla'|u'\equiv 0$. Thus $u'(t,x)=u'(t,x_n)$.
Summarizing the above, we obtain $u(t,x)=u(t,x_n)$, and $u_n(t,x)\equiv 0$.
\end{section}

\begin{section}{The case $\Omega$ is an exterior domain}
 Let $u$ be a bounded very weak ancient solution to (\ref{eq:main1})(\ref{eq:main2}) in an exterior domain $\Omega$ (i.e, the complement of $\Omega$ is homeomorphic to a ball), we show $u(x,t)=f(t)+O(\frac{1}{|x|^{n-2}})$ with some bounded $f$ and $n\ge 3$, in this section. More precisely we have the following theorem:
\begin{theorem}
 Let $u$ be a bounded very weak ancient solution to equations (\ref{eq:main1})(\ref{eq:main2}) in $\Omega\times(-\infty,0)$, where $\Omega\subset R^n$ ($n\ge 3$) is an exterior domain with $\Omega^c\subset B_R$ for some $R>0$. $\|u\|_{L^{\infty}_{x,t}}\leq 1$. Then there exists a function $a(t)$ with $|a(t)|\leq 1$ for almost every $t$ such that
\begin{equation}
 |u(x,t)-a(t)|\leq \frac{C(n,R)}{|x|^{n-2}} \quad\mbox{for almost every $|x|\ge 4R$ and $t<0$.}
\end{equation}
\end{theorem}

For such purpose, we first mollify $u$ in $t$ variable as in lemma \ref{lm:mollification}, it is clear that $u^{\epsilon}$ thus obtained still satisfies equations (\ref{eq:main3})(\ref{eq:main4}). Our first goal is to show that $u^{\epsilon}$ is smooth in $\overline{\Omega}\times (-\infty,0)$ and $u^{\epsilon}|_{\partial\Omega,t<0}=0$.
\begin{lemma}
 Let $u$ and $u^{\epsilon}$ be as above. Then $u^{\epsilon}$ is smooth in $\overline{\Omega}\times (-\infty,0)$ and $u^{\epsilon}|_{\partial\Omega,t<0}=0$.
\end{lemma}

\noindent
\textbf{Proof}: Clearly $u^{\epsilon}$ verifies
\begin{equation}
 \left.\begin{array}{rl}
        \partial_tu^{\epsilon}-\Delta u^{\epsilon}+\nabla q&=0\\
  \mbox{div}~~u^{\epsilon}&=0
       \end{array}\right\}\quad\mbox{in $\Omega\times(-\infty,0)$.}
\end{equation}

In equations (\ref{eq:main3})(\ref{eq:main4}), take test functions as $\eta(t)\phi(x)$, $\eta(t)\psi(x)$ respectively for smooth $\phi$, $\psi$ with ${\rm supp}~~\phi$, ${\rm supp}~~\psi\Subset \overline{\Omega}$, $\mbox{div}~~\phi=0$, $\phi|_{\partial\Omega}=0$ and $\eta\in C_c^{\infty}(-\infty,0)$. We obtain by integration by parts (and definition of $u^{\epsilon}$):
\begin{eqnarray*}
 \int_{-\infty}^0\int_{\Omega}(-\partial_tu^{\epsilon}\phi+u^{\epsilon}\Delta\phi)\eta(t)dxdt&=&0\mbox{,}\\
\int_{-\infty}^0\int_{\Omega}u^{\epsilon}\cdot\nabla\psi\eta(t)dxdt&=&0\mbox{.}
\end{eqnarray*}
Since $\eta$ is arbitrary, we get for any $t\in (-\infty,0)$,
\begin{eqnarray}
 \int_{\Omega}-\partial_tu^{\epsilon}\phi+u^{\epsilon}\Delta\phi dxdt &=&0\mbox{,}\label{eq:no3}\\
\int_{\Omega}u^{\epsilon}\cdot\nabla\psi dxdt &=&0\mbox{.}\label{eq:no4}
\end{eqnarray}
Take $R>0$ sufficiently large such that $\Omega^{c}\subset B_R(0)$. For fixed $t<0$, we can find $v\in C^{1,1/2}(\Omega\cap B_R)$ satisfying
\begin{eqnarray}
 &&\left.\begin{array}{rl}
        -\Delta v+\nabla p&=-\partial_tu^{\epsilon}(\cdot,t)\label{eq:no5}\\
  \mbox{div}~~v&=0
       \end{array}\right\}\quad \mbox{in $\Omega\cap B_R$,}\\
   && v|_{\partial \Omega}=0, \quad v|_{\partial B_R}=u^{\epsilon}(\cdot,t)|_{\partial B_R}\mbox{.}
\end{eqnarray}
Note in the interior of $\Omega$, $u^{\epsilon}$ is smooth by lemma \ref{lm:linearregularity} and definition of $u^{\epsilon}$. The existence of $v$ follows from well-known results of steady Stokes system, we only remark here that the usual no outflow condition required by existence theory is satisfied in our situation and can be easily seen by setting $\psi$ to be 1 in a neighborhood of $\partial \Omega$ in equation (\ref{eq:no4}). Set $w=u^{\epsilon}(\cdot,t)-v$, we claim $w\equiv 0$ in $\Omega\cap B_R$. To prove the claim, take any $\phi\in C^{\infty}(\overline{\Omega\cap B_R})$ with $\mbox{div}~~\phi=0$ and $\phi|_{\partial (\Omega\cap B_R)}=0$, $\psi\in C^{\infty}(\overline{\Omega\cap B_R})$, we write $\phi=\phi_1+\phi_2$, $\psi=\psi_1+\psi_2$ with the following properties: $\mbox{div}~~\phi_1=\mbox{div}~~\phi_2=0$, $\phi_1,~~\phi_2$, $\psi_1,~~\psi_2$ are smooth; $\phi_1$, $\psi_1$ equal $\phi$ and $\psi$ in a neighborhood of $\partial \Omega$, vanishes in a neighborhood of $\partial B_R$ respectively. The existence of such decomposition of $\psi$ is clear. To obtain such this decomposition for $\phi$, one can localize $\phi$ by a standard cutoff function vanishing in a neighborhood of $\partial B_R$, then use Bogovski's theorem to deal with the divergence free condition, we omit the details here. With these decompositions, equations (\ref{eq:no3})(\ref{eq:no4}), the definitions of $v$ and the fact that $u^{\epsilon}$ is smooth away from $\partial \Omega$, we easily obtain: $\int_{\Omega}w\Delta\phi dx=0$ and $\int_{\Omega}w\nabla \psi dx=0$. Thus by result in section 3, this implies $w=0$. Therefore, $u^{\epsilon}(\cdot,t)\in C^{1,1/2}(\overline{B_R\cap \Omega})$ and $u^{\epsilon}|_{\partial\Omega}=0$. A simple boostraping argument gives smoothness of $u^{\epsilon}$. The lemma is proved.\\

\noindent
\textbf{Proof of main result of this section}\\
Let us first summarize the above results as follows:\\
$u^{\epsilon}$ is in $C^{\infty}(\overline{\Omega}\times(-\infty,0))$ with all derivatives bounded (with bounds depending on $\epsilon$) and, $u^{\epsilon}$ satisfies
\begin{eqnarray}
 \left.\begin{array}{rl}
        \partial_tu^{\epsilon}-\Delta u^{\epsilon}+\nabla q&=0\\
 \mbox{div}~~u^{\epsilon}&=0
       \end{array}\right\}&&\quad\mbox{in $\Omega\times(-\infty,0),$}\\
   u^{\epsilon}(\cdot,t)|_{\partial \Omega}=0\mbox{.}&&
\end{eqnarray}
We extend $u^{\epsilon}$ to $R^n$ by setting $u^{\epsilon}=0$ in $\Omega^c$. It is not hard to see the extended $u^{\epsilon}$ satisfies
\begin{equation}
 \left.\begin{array}{rl}
        \partial_tu^{\epsilon}-\Delta u^{\epsilon}+\nabla q&=\mu\\
  \mbox{div}~~u&=0
       \end{array}\right\}\quad\mbox{in $R^n\times(-\infty,0)$.}
  \end{equation}
for the measure $\mu=f^{\epsilon}(x,t)d\sigma$, where $f^{\epsilon}=\frac{\partial u^{\epsilon}}{\partial n}-qn$ on $\partial \Omega$ and $d\sigma$ is the surface measure of $\partial\Omega$. We set 
$$
v^{\epsilon}(x,t)=\int_{-\infty}^tPe^{\Delta (t-s)}\mu(\cdot,s)ds=\int_{-\infty}^t\frac{1}{(t-s)^{n/2}}\int_{\partial \Omega}k(\frac{x-y}{\sqrt{t-s}})f^{\epsilon}(y,s)d\sigma(y)ds\,\,,
$$
where $P$ is the Helmholtz projection to divergence free vector field and $k(\cdot)$ is the kernel of $Pe^{\Delta}$. Thus $|k(y)|\leq \frac{C(n)}{(1+|y|)^n}$. Simple calculations show $v^{\epsilon}$ verifies the following estimates:
 $\|v^{\epsilon}(\cdot,t)\|_{L^1(B_{2R})}\leq C(R,n,\epsilon), \quad \mbox{and}~~ |v^{\epsilon}(x,t)|\leq \frac{C(n,\epsilon,R)}{|x|^{n-2}}$, $|\nabla v^{\epsilon}(x,t)|\leq \frac{C(n,\epsilon,R)}{|x|^{n-1}}$, for $|x|\ge 2R$. In the above calculations $n\ge3$ is important. One can check if $n=2$ the integral might diverge. Clearly $w^{\epsilon}:=u^{\epsilon}-v^{\epsilon}$ satisfies
\begin{equation}
 \left.\begin{array}{rl}
        \partial_tw^{\epsilon}-\Delta w^{\epsilon}+\nabla q&=0\\
  \mbox{div}~~w^{\epsilon}&=0
       \end{array}\right\}\quad\mbox{in $R^n\times(-\infty,0)$.}
  \end{equation}
Thus $w^{\epsilon}=a^{\epsilon}(t)$ and consequently $u^{\epsilon}=v^{\epsilon}+a^{\epsilon}(t)$. At this stage, we would like to pass $\epsilon$ to zero. The decay estimate for $v^{\epsilon}$, however, depends on $\epsilon$ (since the bounds of $f^{\epsilon}$ depends on $\epsilon$). Thus we must first remove this dependence. To do this, let us consider vorticity $\omega^{\epsilon}_{ij}=\partial_iu^{\epsilon}_j-\partial_ju^{\epsilon}_i$ for $1\leq i,j\leq n$. $\omega^{\epsilon}_{ij}$ satisfy $\partial_t\omega^{\epsilon}_{ij}-\Delta \omega^{\epsilon}_{ij}=0$ in $(R^n\backslash B_{2R})\times (-\infty,0)$. By interior regularity of solution to heat equation, the scaling invariance $\omega^{\epsilon}_{ij}(x,t)\to \omega^{\epsilon}_{ij}(Mx,M^2x)$, and the $L^{\infty}$ bound on $u$, we easily conclude:\\

\noindent
\textit{$\omega^{\epsilon}_{ij}$ is smooth in $(R^n\backslash B_{3R})\times (-\infty,0)$ with $|\partial_t\omega^{\epsilon}_{ij}(x,t)|\leq \frac{C(n,R)}{|x|^3}$, $|x|\ge 3R$, the estimate being independent of $\epsilon$.}\\

\noindent
Now fix $\epsilon=\epsilon_0>0$. From estimates of $u^{\epsilon_0}$ above, we know $|\omega^{\epsilon_0}_{ij}(x,t)|\leq \frac{C(n,\epsilon_0,R)}{|x|^{n-1}}$, $|x|\ge 3R$. By estimates of $\partial_t \omega_{ij}$ and definition of $\omega^{\epsilon}_{ij}$, we see
$|\omega_{ij}(x,t)-\omega^{\epsilon_0}_{ij}(x,t)|\leq \frac{C(n,\epsilon_0,R)}{|x|^3}$, $|x|\ge 3R$. Therefore, $|\omega_{ij}(x,t)|\leq C(n,\epsilon_0,R)(\frac{1}{|x|^3}+\frac{1}{|x|^{n-1}})$, $|x|\ge 3R$. In fact, by a boostrap argument (better decay estimate of $\omega_{ij}$ improves estimate of $\partial_t\omega_{ij}$ which in turn implies better decay estimate of $\omega_{ij}$), one can upgrade this estimate to $|\omega_{ij}(x,t)|\leq \frac{C(n,\epsilon_0,R)}{|x|^{n-1}}$. Thus $|\omega^{\epsilon}(x,t)|\leq \frac{C(n,\epsilon_0,R)}{|x|^{n-1}}$, (now independent of $\epsilon$) for $|x|\ge 3R$. For each fixed $t\in (-\infty,0)$, $u^{\epsilon}_i$ solves
\begin{equation}
 -\Delta u^{\epsilon}_i=-\sum_{j=1}^n\partial_j\omega^{\epsilon}_{ij} \quad\mbox{in $R^n\backslash B_{3R}$.}
\end{equation}
From this and the boundedness of $u$, it is not hard to see $u^{\epsilon}_i$ verifies the following bound: $|u^{\epsilon}(x,t)-a^{\epsilon}(t)|\leq \frac{C(n,\epsilon_0,R)}{|x|^{n-2}}$ for $|x|\ge 4R$. Passing $\epsilon\to 0$ the conclusion of this section is reached.

\end{section}

\end{document}